\documentclass[reqno,12pt]{amsart}
\usepackage[T1]{fontenc}
\usepackage[utf8]{inputenc}
\usepackage[english]{babel}
\usepackage{amssymb,amsmath,amsthm,amsfonts,xcolor,enumerate,hyperref,comment,longtable,cleveref,mathtools}

\usepackage{times}
\usepackage{cite}
\usepackage{pdflscape}
\usepackage[mathcal]{euscript}
\usepackage{tikz}
\usepackage{cancel}
\usepackage{stmaryrd}
\usepackage{longtable}

\usepackage[a4paper,top=3cm, bottom=3cm, left=3cm, right=3cm]{geometry}

\theoremstyle{plain}
\newtheorem{theorem}[subsection]{Theorem}
\newtheorem{conjecture}[subsection]{Conjecture}

\theoremstyle{definition}
\newtheorem{definition}[subsection]{Definition}
\theoremstyle{remark}

\DeclareMathOperator{\gr}{gr}

\begin{document}
\sloppy
	
	\title[Local and 2-local automorphisms of some solvable Leibniz algebras]{Local and 2-local automorphisms of some solvable Leibniz algebras}

	\author{F.N.Arzikulov, I.A.Karimjanov}
\address[F.N.Arzikulov, I.A.Karimjanov]	{Institute of Mathematics Uzbekistan Academy of Sciences, 4, University street, Olmazor, Tashkent, 100174, Uzbekistan
\newline
		and\newline
	Department of Mathematics, Andijan State University, 129, Universitet Street, Andijan, 170100, Uzbekistan}
	\email{arzikulovfn@rambler.ru, iqboli@gmail.com}

	\author{S.M.Umrzaqov}
	\address[S.M.Umrzaqov]{	
		Department of Mathematics, Andijan State University, 129, Universitet Street, Andijan,
170100, Uzbekistan }
	\email{sardor.umrzaqov1986@gmail.com}

	
	\begin{abstract}
In this paper we prove that any local automorphism on the solvable Leibniz algebras with null-filiform and naturally graded non-Lie filiform nilradicals, whose dimension of complementary space is maximal is an automorphism. Furthermore, the same problem concerning 2-local automorphisms of such algebras is investigated and we obtain the analogously results for 2-local automorphisms.
	\end{abstract}
	
	
	\maketitle
	
	\section{Introduction}

 The history of local mappings begins
with the Gleason-Kahane-\.{Z}elazko theorem in \cite{AMG} and \cite{JPK_WZ}, which is a fundamental contribution in the theory of Banach algebras.
This theorem asserts that every unital linear functional $F$ on a complex unital Banach algebra $A$,such that $F(a)$ belongs to the spectrum $\sigma(a)$
of $a$ for every $a\in A,$ is multiplicative. In modern terminology this is equivalent to the following condition: every unital linear local homomorphism from a unital complex Banach algebra $A$ into ${\Bbb C}$ is multiplicative. We recall that a linear map $T$ from a Banach algebra $A$ into a Banach algebra $B$ is said to be a local homomorphism if for every $a$ in $A$ there exists a homomorphism $\Phi_a : A\to B$, depending on $a$,
such that $T(a)=\Phi_a(a)$.

Later, in \cite{Kad}, R. Kadison introduces the concept of local derivation and proves
that each continuous local derivation from a von Neumann algebra into its dual Banach
bemodule is a derivation. B. Jonson \cite{Jon} extends the above result by proving that every
local derivation from a C*-algebra into its Banach bimodule is a derivation. In particular, Johnson
gives an automatic continuity result by proving that local derivations of a C*-algebra $A$ into a
Banach $A$-bimodule $X$ are continuous even if not assumed a priori to be so
(cf. \cite[Theorem 7.5]{Jon}). Based on these results, many authors have studied
local derivations on operator algebras.

A similar notion, which characterizes non-linear generalizations of automorphisms, was introduced
by  \v{S}emrl in \cite{S} as $2$-local automorphisms.
He described such maps on the algebra $B(H)$ of all bounded linear operators on an infinite dimensional separable Hilbert space $H$.
After the work of \v{S}emrl, it is appeared numerous new results related to the description of local and $2$-local automorphisms of algebras
(see, for example, \cite{AyupovKudaybergenov}, \cite{akoz}, \cite{ChenWang}, \cite{Costantini}, \cite{KK}, \cite{AKK}).

Leibniz algebra is a generalization of Lie algebra in natural way. Leibniz algebras have been defined by Loday in \cite{lod} as a non-antisymmetric version of Lie algebras. The problem of classification of finite-dimensional Leibniz algebras is fundamental and a very complicated problem.	Last 30 years the Leibniz algebras has been actively investigated and a lot of papers have been devoted to the study of these algebras \cite{omir1, ayupov, Sib, kar}. The analogue of the Levi-Malcev decomposition for Leibniz algebras  was proved by D.W. Barnes \cite{barn}, that asserts that any Leibniz algebra decomposes into a semidirect sum of its solvable radical and a semisimple Lie algebra.
The semisimple part can be described from the simple Lie ideals, therefore,  the main problem of the description of the finite-dimensional Leibniz algebras consists of the study of the solvable Leibniz algebras. Then, a lot of progress has been made in the study of classifications concerning solvable Leibniz algebras with a given nilradicals \cite{NF,F2,kam}.

In the paper \cite{AyupovKudaybergenov} the authors proved that every 2-local automorphism on a finite-dimensional semi-simple Lie algebra  over an algebraically closed field of characteristic zero is an automorphism and showed that each finite-dimensional nilpotent Lie algebra  with  dimension $\geq2$ admits a 2-local automorphism which is not an automorphism. Later by Ayupov, Kudaybergenov and Omirov proved that every 2-local automorphism on a complex finite-dimensional simple Leibniz algebra is an automorphism and show that nilpotent Leibniz algebras admit 2-local automorphisms which are not automorphisms. A similar problem concerning local automorphism on simple Leibniz algebras is reduced to the case of simple Lie algebras \cite{AKK}. The present paper is devoted to local automorphisms on solvable Leibniz algebras.

This paper is organized as follows. In Sect. \ref{S:pre}, we provide some basic
concepts needed for this study. In Sect. \ref{S:big}, we investigate local automorphisms on solvable Leibniz algebras with null-filiform and naturally graded non-Lie filiform nilradicals. The last section is devoted to
2-local automorphisms on such type solvable Leibniz algebras. Finally, we give conjecture that the local and 2-local automorphisms on the solvable Leibniz algebras with a given nilradical the dimension of whose complementary space is maximal is an automorphism.
	
	\section{Preliminaries}\label{S:pre}
	
\begin{definition} An algebra $L$ over a field $\mathbb{K}$ is called a Leibniz algebra if for any $x,y,z\in L$, the Leibniz identity
		\[[[x,y],z]=[[x,z],y]+[x,[y,z]]\] is satisfied, where $[-,-]$ is the multiplication in $L$.
	\end{definition}

For a Leibniz algebra $L$ we consider the following derived and lower central series:
	\begin{align*}
		& \text{(i)}   &L^{(1)}= & \ L, \  & L^{(n+1)}= & \ [L^{(n)},L^{(n)}],  & n>1; \\
		& \text{(ii)}  & L^1= &  \ L, \ & L^{n+1}=&  \ [L^n,L],  & n>1.
	\end{align*}

	\begin{definition} An algebra $L$ is called solvable (nilpotent) if there exists $s\in \mathbb{N}$ \ ($k\in \mathbb{N}$, respectively) such that $L^{(s)}=0$ \ ($L^k=0$, respectively).
The minimal number $s$ (respectively, $k$) with such property is called index of
solvability (respectively, of nilpotency) of the algebra $L.$
	\end{definition}

Evidently, the index of nilpotency of  an $n$-dimensional algebra is not greater than $n+1.$

\begin{definition} An $n$-dimensional Leibniz algebra is called null-filiform if $\dim L^i=n+1-i, \ 1\leq i \leq n+1.$
\end{definition}

Actually, a nilpotent Leibniz algebra is null-filiform if it is a one-generated algebra. Note, that this notion has no sense in Lie algebras case, because they are at least two-generated.

\begin{definition} A Leibniz algebra $L$ is said to be filiform if
$\dim L^i=n-i$, for $2\leq i \leq n$, where $n=\dim L$.
\end{definition}

\begin{definition} Given a filiform Leibniz algebra $L,$ put
$L_i=L^i/L^{i+1}, \ 1 \leq i\leq n-1,$ and $\gr L = L_1 \oplus
L_2\oplus\dots L_{n-1}.$ Then $[L_i,L_j]\subseteq L_{i+j}$ and we
obtain the graded algebra $\gr L$. If $\gr L$ and $L$ are isomorphic,
denoted by $\gr L\cong L,$ we say that the algebra $L$ is naturally
graded.
\end{definition}

\begin{definition} The  (unique) maximal nilpotent ideal of a Leibniz algebra is called the nilradical of the algebra.
\end{definition}

Let $R$ be a solvable Leibniz algebra. Then it can be decomposed into the form $R=N \oplus Q$, where $N$ is the  nilradical and $Q$ is the complementary vector space. Since the square of a solvable algebra is a nilpotent ideal and the finite sum of nilpotent ideals is a nilpotent ideal too, then the ideal $R^2$ is nilpotent, i.e. $R^2\subseteq N$ and consequently, $Q^2\subseteq N$.

Now, we present the classification results for arbitrary dimensional solvable Leibniz algebras with null-filiform and naturally graded non-Lie filiform nilradicals, whose the dimension of complementary space is maximal.

\begin{theorem} \cite{NF} Let $R_0$ be a solvable Leibniz algebra with null-filiform nilradical. Then there exists a basis $\{e_0, e_1, e_2, \dots, e_n\}$ of the algebra $R_0$ such that the multiplication table of $R_0$ with respect to this basis has the following form:
\[R_0: \left\{ \begin{aligned}
{}[e_i,e_1] & =e_{i+1}, && 0\leq i\leq n-1,\\
[e_i,e_0] & =-ie_i, && 1\leq i\leq n.
\end{aligned}\right.\]
\end{theorem}

\begin{theorem}\cite{kam,F2}\label{F1n} An arbitrary $(n+2)$-dimensional solvable Leibniz algebra with naturally graded non-Lie filiform nilradical is isomorphic to one of the following non-isomorphic algebras:
\[R_1: \left\{\begin{array}{ll}
[e_i,e_1]=e_{i+1}, &  2\leq i \leq {n-1},  \\[1mm]
[e_1,x]=-[x,e_1]=e_1,      &               \\[1mm]
[e_i,x]=(i-1)e_i, &  2\leq i \leq n, \\[1mm]
[e_i,y]=e_i,       &  2\leq i\leq n,     \\[1mm]
\end{array}\right.\]
\[ R_2: \left\{\begin{array}{ll}
[e_1,e_1]=e_3, & \\[1mm]
[e_i,e_1]=e_{i+1}, & 3\leq i\leq n-1,\\[1mm]
[e_1,x]=-[x,e_1]=e_1, & \\[1mm]
[e_i,x]=(i-1)e_i, & 3\leq i\leq n, \\[1mm]
[e_2,y]=-[y,e_2]=e_2, &
\end{array}\right.\]
\[ R_3: \left\{\begin{array}{ll}
[e_1,e_1]=e_3, & \\[1mm]
[e_i,e_1]=e_{i+1}, & 3\leq i\leq n-1,\\[1mm]
[e_1,x]=-[x,e_1]=e_1, & \\[1mm]
[e_i,x]=(i-1)e_i, & 3\leq i\leq n, \\[1mm]
[e_2,y]=e_2 &
\end{array}\right.\]
where $\{e_1, \dots, e_n, x, y\}$ is a basis of the algebra.
\end{theorem}

An automorphism is simply a bijective homomorphism of an object with itself. Now we give the definitions of local and 2-local automorphisms.

\begin{definition}
Let $A$ be an algebra. A linear map $\Phi : A \to A$ is called a local automorphism, if for
any element $x \in A$ there exists an automorphism $\varphi_x : A \to A$ such that $\Phi(x) = \varphi_x(x)$.
\end{definition}

\begin{definition}
	$A$ (not necessary linear) map $\phi : A \rightarrow A $ is called a 2-local automorphism, if
for any elements $x, y \in A$ there exists an automorphism $\phi_{x,y}: A \rightarrow A$ such that
$\phi(x) = A_{x,y}(x), \phi(y) = A_{x,y}(y).$
\end{definition}

Below we give the descriptions of automorphisms on solvable Leibniz algebras $R_0, R_1, R_2$ and $R_3$.

\begin{theorem} \label{1} \cite{F1} A linear map $\varphi:R_0\to R_0$ is an automorphism if and only if $\varphi$ has the
following form:
\[\varphi(e_i)=\sum\limits_{j=i}^n \frac{\alpha^{j-i} \beta^{i}}{(j-i)!} e_j, \quad   0\leq i \leq n, \]
where $\beta\neq 0$.
\end{theorem}

\begin{theorem} \label{11} \cite{F1} A linear maps $\varphi_1, \varphi_2$ and $\varphi_3$ are automorphisms on algebras $R_1, R_2$ and $R_3$ respectively if and only if when $\varphi_1, \varphi_2$ and $\varphi_3$ have the
following forms:
\[\left\{\begin{array}{ll}
     \varphi_1(e_1)=\alpha e_1, & \\
      \varphi_1(e_i)=\sum\limits_{j=i}^n \frac{(-1)^{j-i}\alpha^{i-2} \beta\gamma^{j-i}}{(j-i)!} e_j,&  2\leq i\leq n,\\
       \varphi_1(x)=\gamma e_1+x, & \\
        \varphi_1(y)=y, & \\
     \end{array}\right.\]
where $\alpha\beta\neq 0$,
     \[\left\{\begin{array}{ll}
     \varphi_2(e_1)=\alpha e_1+\sum\limits_{i=3}^n\frac{(-1)^i\alpha\beta^{i-2}}{(i-2)!}e_i, & \\
       \varphi_2(e_2)=\gamma e_2, & \\
      \varphi_2(e_i)=\sum\limits_{j=i}^n \frac{(-1)^{j-i}\alpha^{i-1} \beta^{j-i}}{(j-i)!} e_j,&  3\leq i\leq n,\\
       \varphi_2(x)=\beta e_1+\sum\limits_{i=3}^n\frac{(-1)^i\beta^{i-1}}{(i-1)!}e_i+x, & \\
        \varphi_2(y)=\delta e_2+y & \\
     \end{array}\right.\]
where $\alpha\gamma\neq 0$,
     \[\left\{\begin{array}{ll}
     \varphi_3(e_1)=\alpha e_1+\sum\limits_{i=3}^n\frac{(-1)^i\alpha\beta^{i-2}}{(i-2)!}e_i, & \\
       \varphi_3(e_2)=\gamma e_2, & \\
      \varphi_3(e_i)=\sum\limits_{j=i}^n \frac{(-1)^{j-i}\alpha^{i-1} \beta^{j-i}}{(j-i)!} e_j,&  3\leq i\leq n,\\
       \varphi_3(x)=\beta e_1+\sum\limits_{i=3}^n\frac{(-1)^i\beta^{i-1}}{(i-1)!}e_i+x, & \\
        \varphi_3(y)=y, & \\
     \end{array}\right.\]
where $\alpha\gamma\neq 0$.
\end{theorem}

	\section{Local automorhisms of Solvable Leibniz algebras}\label{S:big}
	
\begin{theorem}
Every local automorphism of $R_0$ is an automorphism.
\end{theorem}

\begin{proof}
Let $\Phi$ be an arbitrary local automorphism of $R_0$.
By the definition for all $x\in R_0$ there exists an automorphism $\varphi_x$ on $R_0$ such that
$$
\Phi(x)=\varphi_x(x).
$$
By theorem \ref{1},  the automorphism $\varphi_x$ has the following matrix form:
\[A_x=\begin{pmatrix}
        1&0&0&\dots&0&0\\
        \alpha_x&\beta_x&0&\dots&0&0\\
        \frac{\alpha_x^2}{2}&\alpha_x\beta_x&\beta_x^2&\dots&0&0\\
        \vdots&\vdots&\vdots&\ddots&\vdots&\vdots\\
         \frac{\alpha_x^{n-1}}{(n-1)!}&\frac{\alpha_x^{n-2}\beta_x}{(n-2)!}&\frac{\alpha_x^{n-3}\beta_x^2}{(n-3)!}&\dots&\beta_x^{n-1}&0\\
         \frac{\alpha_x^n}{n!}&\frac{\alpha_x^{n-1}\beta_x}{(n-1)!}&\frac{\alpha_x^{n-2}\beta_x^2}{(n-2)!}&\dots&\alpha_x\beta_x^{n-1}&\beta_x^n
\end{pmatrix}.\]

Let $A$ be the matrix of $\Phi$
then by choosing subsequently $x=e_0, x=e_1, \ldots, x=e_n$ and using $\Phi(x)=\varphi_x(x),$ i.e. $A\bar{x}=A_x\bar{x},$ where $\bar{x}$ is the vector
corresponding to $x$ and
\[A=\begin{pmatrix}
        1&0&0&\dots&0&0\\
        \alpha&\beta&0&\dots&0&0\\
        \frac{\alpha^2}{2}&\alpha\beta&\beta^2&\dots&0&0\\
        \vdots&\vdots&\vdots&\ddots&\vdots&\vdots\\
         \frac{\alpha^{n-1}}{(n-1)!}&\frac{\alpha^{n-2}\beta}{(n-2)!}&\frac{\alpha^{n-3}\beta^2}{(n-3)!}&\dots&\beta^{n-1}&0\\
         \frac{\alpha^n}{n!}&\frac{\alpha^{n-1}\beta}{(n-1)!}&\frac{\alpha^{n-2}\beta^2}{(n-2)!}&\dots&\alpha\beta^{n-1}&\beta^n
\end{pmatrix},\]
it is easy to see that
\[A=\begin{pmatrix}
        1&0&0&\dots&0&0\\
        \alpha_{e_0}&\beta_{e_1}&0&\dots&0&0\\
        \frac{\alpha_{e_0}^2}{2}&\alpha_{e_1}\beta_{e_1}&\beta_{e_2}^2&\dots&0&0\\
        \vdots&\vdots&\vdots&\ddots&\vdots&\vdots\\
         \frac{\alpha_{e_0}^{n-1}}{(n-1)!}&\frac{\alpha_{e_1}^{n-2}\beta_{e_1}}{(n-2)!}&\frac{\alpha_{e_2}^{n-3}\beta_{e_2}^2}{(n-3)!}&\dots&\beta_{e_{n-1}}^{n-1}&0\\
         \frac{\alpha_{e_0}^n}{n!}&\frac{\alpha_{e_1}^{n-1}\beta_{e_1}}{(n-1)!}&\frac{\alpha_{e_2}^{n-2}\beta_{e_2}^2}{(n-2)!}&\dots&\alpha_{e_{n-1}}\beta_{e_{n-1}}^{n-1}&\beta_{e_n}^n
\end{pmatrix}.\]

Since $\Phi$ is linear we have
\[
\Phi(x+y)=\Phi(x)+\Phi(y), \quad \forall x,y\in R_0.   \eqno{(3.1)}
\]

Consider the equality
\[\Phi(e_0+e_k)=\sum\limits_{j=0}^n\frac{\alpha_{e_0+e_k}^j}{j!}e_j+\sum\limits_{j=k}^n\frac{\alpha_{e_0+e_k}^{j-k}\beta_{e_0+e_k}^k}{(j-k)!}e_j, \quad 2\leq k\leq n-1.\]

On the other hand, we have
\[\Phi(e_0+e_k)=\Phi(e_0)+\Phi(e_k)=\sum\limits_{j=0}^n\frac{\alpha_{e_0}^j}{j!}e_j+\sum\limits_{j=k}^n\frac{\alpha_{e_k}^{j-k}\beta_{e_k}^k}{(j-k)!}e_j, \quad 2\leq k\leq n-1.\]

Comparing coefficients of the basis elements, we derive:
\[\alpha_{e_0+e_k}=\alpha_{e_0}, \quad \beta_{e_0+e_k}=\beta_{e_k}, \quad \alpha_{e_0+e_k}=\alpha_{e_k}, \quad 2\leq k\leq n-1.\]

Which implies \[\alpha_{e_0}=\alpha_{e_k}, \quad 2\leq k\leq n-1.\]

From equality (3.1), we have
\[\Phi(e_1+e_k)=\sum\limits_{j=1}^n\frac{\alpha_{e_1+e_k}^{j-1}\beta_{e_1+e_k}}{(j-1)!}e_j+\sum\limits_{j=k}^n\frac{\alpha_{e_1+e_k}^{j-k}\beta_{e_1+e_k}^k}{(j-k)!}e_j, \quad 3\leq k\leq n-1.\]

On the other hand, we obtain
\[\Phi(e_1+e_k)=\Phi(e_1)+\Phi(e_k)=\sum\limits_{j=1}^n\frac{\alpha_{e_1}^{j-1}\beta_{e_1}}{j!}e_j+\sum\limits_{j=k}^n\frac{\alpha_{e_k}^{j-k}\beta_{e_k}^k}{(j-k)!}e_j, \quad 3\leq k\leq n-1.\]

From the previous equalities, we deduce:
\[\beta_{e_1+e_k}=\beta_{e_1}, \quad \alpha_{e_1+e_k}=\alpha_{e_1}, \quad \beta_{e_1+e_k}=\beta_{e_k}, \quad \alpha_{e_1+e_k}=\alpha_{e_k}, \quad 3\leq k\leq n-1, \quad \text{i.e.}\]
\[\alpha_{e_1}=\alpha_{e_k}, \quad \beta_{e_1}=\beta_{e_k}, \qquad 3\leq k\leq n-1.\]

With a similar argument, we obtain
\[\Phi(e_2+e_k)=\sum\limits_{j=2}^n\frac{\alpha_{e_2+e_k}^{j-2}\beta_{e_2+e_k}}{(j-2)!}e_j+\sum\limits_{j=k}^n\frac{\alpha_{e_2+e_k}^{j-k}\beta_{e_2+e_k}^k}{(j-k)!}e_j, \quad 4\leq k\leq n-1.\]
and
\[\Phi(e_2+e_k)=\Phi(e_2)+\Phi(e_k)=\sum\limits_{j=2}^n\frac{\alpha_{e_2}^{j-2}\beta_{e_2}}{(j-2)!}e_j+\sum\limits_{j=k}^n\frac{\alpha_{e_k}^{j-k}\beta_{e_k}^k}{(j-k)!}e_j, \quad 4\leq k\leq n-1.\]
and hence
\[\beta_{e_2+e_k}=\beta_{e_2}, \quad \alpha_{e_2+e_k}=\alpha_{e_2}, \quad \beta_{e_2+e_k}=\beta_{e_k}, \quad \alpha_{e_2+e_k}=\alpha_{e_k}, \quad 4\leq k\leq n-1, \quad \text{i.e.}\]
\[\alpha_{e_2}=\alpha_{e_k}, \quad \beta_{e_2}=\beta_{e_k}, \qquad 4\leq k\leq n-1.\]

Finally, from
\[\Phi(e_1+e_n)=\sum\limits_{j=1}^n\frac{\alpha_{e_1+e_n}^{j-1}\beta_{e_1+e_n}}{(j-1)!}e_j+\beta_{e_1+e_n}^ne_n, \]
and
\[\Phi(e_1+e_n)=\Phi(e_1)+\Phi(e_n)=\sum\limits_{j=1}^n\frac{\alpha_{e_1}^{j-1}\beta_{e_1}}{(j-1)!}e_j+\beta_{e_n}^ne_n,\]
we obtain $\beta_{e_1}=\beta_{e_n}$.

Thus, we obtain that the local automorphism $\Phi$ has the following form:
\[\Phi(e_i)=\sum\limits_{j=i}^n \frac{\alpha^{j-i}_{e_0} \beta^{i}_{e_1}}{(j-i)!} e_j, \quad   0\leq i \leq n. \]
Note that, by the definition of a local automorphism, $\beta_{e_1}\neq 0$.
Hence, by theorem \ref{1}, $\Phi$ is an automorphism. This ends the proof.
\end{proof}

\begin{theorem}
Every local automorphism of $R_1$ is an automorphism.
\end{theorem}
\begin{proof}
 By applying the similar arguments used above we can assume the local automorphism $\Phi$ on $R_1$ has the following matrix:
\[\begin{pmatrix}
\alpha_{e_1}&               0                       &                       0                          &           \cdots & 0 & \gamma_{x} & 0 \\
    0       &        \beta_{e_2}                    &                       0                          &          \cdots & 0 &    0   & 0 \\
    0       &     -\beta_{e_2}\gamma_{e_2}          &           \alpha_{e_3}\beta_{e_3}                &           \cdots & 0 & 0 & 0 \\
       \vdots       &               \vdots                       &                       \vdots          & \ddots &  \vdots  & \vdots  & \vdots \\
    0       & \frac{(-1)^{n-2}\beta_{e_2}\gamma_{e_2}^{n-2}}{(n-2)!} &\frac{(-1)^{n-3}\alpha_{e_3}\beta_{e_3}\gamma_{e_3}^{n-3}}{(n-3)!}  & \cdots & \alpha_{e_n}^{n-2}\beta_{e_n} &  0   &  0   \\
    0       &               0                       &                       0                          &          \cdots & 0 & 1 & 0 \\
    0       &               0                       &                       0                          &         \cdots & 0 & 0 & 1 \\
  \end{pmatrix}.\]

Since $\Phi$ is linear, we have
\[\sum\limits_{j=i}^n \frac{(-1)^{j-i}\alpha_{e_i+x}^{i-2} \beta_{e_i+x}\gamma_{e_i+x}^{j-i}}{(j-i)!} e_j+\gamma_{e_i+x}e_1+x=\]
\[\Phi(e_i+x)=\Phi(e_i)+\Phi(x)=\sum\limits_{j=i}^n \frac{(-1)^{j-i}\alpha_{e_i}^{i-2} \beta_{e_i}\gamma_{e_i}^{j-i}}{(j-i)!} e_j+\gamma_{x}e_1+x,\]
where $ 2\leq i\leq n-1.$

Follows, we obtain
\[\alpha_{e_i+x}^{i-2} \beta_{e_i+x}=\alpha_{e_i}^{i-2} \beta_{e_i}, \quad \alpha_{e_i+x}^{i-2} \beta_{e_i+x}\gamma_{e_i+x}=\alpha_{e_i}^{i-2} \beta_{e_i}\gamma_{e_i}, \quad \gamma_{e_i+x}=\gamma_{x}.\]
Hence,
\[\gamma_{e_i+x}=\gamma_x, \quad \gamma_{e_i+x}=\gamma_{e_i}\]
and
\[\gamma_x=\gamma_{e_i}, \quad 2\leq i\leq n-1.\]

With a similar argument
\[\alpha_{e_1+e_2+e_i}e_1+\sum\limits_{j=2}^n \frac{(-1)^{j-2} \beta_{e_1+e_2+e_i}\gamma_{e_1+e_2+e_i}^{j-2}}{(j-2)!} e_j+\sum\limits_{j=i}^n \frac{(-1)^{j-i}\alpha_{e_1+e_2+e_i}^{i-2} \beta_{e_1+e_2+e_i}\gamma_{e_1+e_2+e_i}^{j-i}}{(j-i)!} e_j=\]
\[\Phi(e_1+e_2+e_i)=\Phi(e_1)+\Phi(e_2)+\Phi(e_i)=\]\[\alpha_{e_1}e_1+\sum\limits_{j=2}^n \frac{(-1)^{j-2} \beta_{e_2}\gamma_{e_2}^{j-2}}{(j-2)!} e_j+\sum\limits_{j=i}^n \frac{(-1)^{j-i}\alpha_{e_i}^{i-2} \beta_{e_i}\gamma_{e_i}^{j-i}}{(j-i)!} e_j,\]
we have
\[
\alpha_{e_1+e_2+e_i}=\alpha_{e_1}, \beta_{e_1+e_2+e_i}=\beta_{e_2}, \gamma_{e_1+e_2+e_i}=\gamma_{e_2},   \eqno{(3.2)}
\]
\[
\alpha_{e_1+e_2+e_i}^{i-2} \beta_{e_1+e_2+e_i}=\alpha_{e_i}^{i-2} \beta_{e_i}
\]
which implies
\[\alpha_{e_i}^{i-2} \beta_{e_i}=\alpha_{e_1}^{i-2} \beta_{e_2}, \quad 4\leq i\leq n.\]

Finally, from
\[\begin{array}{c}
\alpha_{e_1+e_3+e_5}e_1+\sum\limits_{j=3}^n \frac{(-1)^{j-3} \alpha_{e_1+e_3+e_5}\beta_{e_1+e_3+e_5}\gamma_{e_1+e_3+e_5}^{j-3}}{(j-3)!} e_j+\\
+\sum\limits_{j=5}^n \frac{(-1)^{j-5}\alpha_{e_1+e_3+e_5}^3 \beta_{e_1+e_3+e_5}\gamma_{e_1+e_3+e_5}^{j-5}}{(j-5)!} e_j=\Phi(e_1+e_3+e_5)=\\
\Phi(e_1)+\Phi(e_3)+\Phi(e_5)=\alpha_{e_1}e_1+\sum\limits_{j=3}^n \frac{(-1)^{j-3}\alpha_{e_3} \beta_{e_3}\gamma_{e_3}^{j-3}}{(j-3)!} e_j+\sum\limits_{j=5}^n \frac{(-1)^{j-5}\alpha_{e_5}^3 \beta_{e_5}\gamma_{e_5}^{j-5}}{(j-5)!} e_j
\end{array}\]
it follows that
\[
\alpha_{e_1+e_3+e_5}=\alpha_{e_1},\quad \alpha_{e_1+e_3+e_5}\beta_{e_1+e_3+e_5}=\alpha_{e_3}\beta_{e_3}, \quad \gamma_{e_1+e_3+e_5}=\gamma_{e_3}.
\]
Using (3.2) for $i=5$ we obtain
\[
\alpha_{e_3}\beta_{e_3}=\alpha_{e_1}\beta_{e_2}.
\]

So, the local automorphism $\Phi$ has the following form:
\[\left\{\begin{array}{ll}
     \Phi(e_1)=\alpha_{e_1} e_1, & \\
      \Phi(e_i)=\sum\limits_{j=i}^n \frac{(-1)^{j-i}\alpha_{e_1}^{i-2} \beta_{e_2}\gamma_x^{j-i}}{(j-i)!} e_j,&  2\leq i\leq n,\\
       \Phi(x)=\gamma_x e_1+x, & \\
        \Phi(y)=y. & \\
     \end{array}\right.\]
By the definition of a local automorphism $\alpha_{e_1}\neq 0$ and $\beta_{e_2}\neq 0$.
Therefore, from theorem \ref{11} we obtain that $\Phi$ is an automorphism.
\end{proof}

\begin{theorem} \label{12}
Every local automorphism of $R_2$ is an automorphism.
\end{theorem}

\begin{proof}
Let $\Phi$ be an arbitrary local automorphism of $R_2$.
By the definition for all $z\in R_2$ there exists an automorphism $\varphi_z$ on $R_2$ such that
$$
\Phi(z)=\varphi_z(z).
$$
By theorem \ref{11} and applying the similar arguments used above we can assume the local automorphism $\Phi$ on $R_2$ has the following matrix:
\[\begin{pmatrix}
        \alpha_{e_1}&0&0&0&\dots&0&\beta_x&0\\
        0&\gamma_{e_2}&0&0&\dots&0&0&\delta_y\\
        -\alpha_{e_1}\beta_{e_1}&0&\alpha_{e_3}^2&0&\dots&0&-\frac{\beta_x^2}{2}&0\\
         \frac{\alpha_{e_1}\beta_{e_1}^2}{2}&0&-\alpha_{e_3}^2\beta_{e_3}&\alpha_{e_4}^3&\dots&0&\frac{\beta_x^3}{6}&0\\
        \vdots&\vdots&\vdots&\vdots&\ddots&\vdots&\vdots&\vdots\\
         \frac{(-1)^n\alpha_{e_1}\beta_{e_1}^{n-2}}{(n-2)!}&0&\frac{(-1)^{n-3}\alpha_{e_3}^2\beta_{e_3}^{n-3}}{(n-3)!}&\frac{(-1)^{n-4}\alpha_{e_4}^3\beta_{e_4}^{n-4}}{(n-4)!}&\dots&\alpha_{e_n}^{n-1}&\frac{(-1)^n\beta_x^{n-1}}{(n-1)!}&0\\
         0&0&0&0&\dots&0&1&0\\
           0&0&0&0&\dots&0&0&1
\end{pmatrix}.\]

Since $\Phi$ is linear we have
\[\alpha_{e_1+e_k}e_1+\sum\limits_{j=3}^n\frac{(-1)^j\alpha_{e_1+e_k}\beta_{e_1+e_k}^{j-2}}{(j-2)!}e_j+\sum\limits_{j=k}^n\frac{(-1)^{j-k}\alpha_{e_1+e_k}^{k-1}\beta_{e_1+e_k}^{j-k}}{(j-k)!}e_j=\]
\[=\Phi(e_1+e_k)=\Phi(e_1)+\Phi(e_k)=\alpha_{e_1}e_1+\sum\limits_{j=3}^n\frac{(-1)^j\alpha_{e_1}\beta_{e_1}^{j-2}}{(j-2)!}e_j+\sum\limits_{j=k}^n\frac{(-1)^{j-k}\alpha_{e_k}^{k-1}\beta_{e_k}^{j-k}}{(j-k)!}e_j\]
for $4\leq k\leq n$.

Comparing coefficients at the basis elements we obtain that
\[\alpha_{e_1+e_k}=\alpha_{e_1}, \quad \beta_{e_1+e_s}=\beta_{e_1}, \quad \alpha_{e_1+e_k}=\alpha_{e_k}, \quad \beta_{e_1+e_s}=\beta_{e_s}, \quad 4\leq k\leq n, \quad 4\leq s\leq n-1.\]

Implies
\[\alpha_{e_1}=\alpha_{e_k}, \quad  \beta_{e_1}=\beta_{e_s},  \quad 4\leq k\leq n, \quad 4\leq s\leq n-1.\]

From the chain of equalities
\[\sum\limits_{j=3}^n\frac{(-1)^{j-3}\alpha_{e_3+e_5}^{2}\beta_{e_3+e_5}^{j-3}}{(j-3)!}e_j+\sum\limits_{j=5}^n\frac{(-1)^{j-5}\alpha_{e_3+e_5}^{4}\beta_{e_3+e_5}^{j-5}}{(j-5)!}e_j=\]
\[=\Phi(e_3+e_5)=\Phi(e_3)+\Phi(e_5)=\sum\limits_{j=3}^n\frac{(-1)^{j-3}\alpha_{e_3}^{2}\beta_{e_3}^{j-3}}{(j-3)!}e_j+\sum\limits_{j=5}^n\frac{(-1)^{j-5}\alpha_{e_5}^{4}\beta_{e_5}^{j-5}}{(j-5)!}e_j.\]

From the previous equalities we deduce that
\[\alpha_{e_3+e_5}=\alpha_{e_3}, \quad \beta_{e_3+e_5}=\beta_{e_3},\quad \alpha_{e_3+e_5}=\alpha_{e_5}, \quad  \beta_{e_3+e_5}=\beta_{e_5}, \text{i.e.}\]
\[\alpha_{e_3}=\alpha_{e_5}, \quad \beta_{e_3}=\beta_{e_5} .\]

Similarly, from
\[\beta_{x+e_k} e_1+\sum\limits_{i=3}^n\frac{(-1)^i\beta_{x+e_k}^{i-1}}{(i-1)!}e_i+x+\sum\limits_{j=k}^n\frac{(-1)^{j-k}\alpha_{x+e_k}^{k-1}\beta_{x+e_k}^{j-k}}{(j-k)!}e_j=\]
\[\Phi(x+e_k)=\Phi(x)+\Phi(e_k)=\beta_{x} e_1+\sum\limits_{i=3}^n\frac{(-1)^i\beta_{x}^{i-1}}{(i-1)!}e_i+x+\sum\limits_{j=k}^n\frac{(-1)^{j-k}\alpha_{e_k}^{k-1}\beta_{e_k}^{j-k}}{(j-k)!}e_j, \]
where $4\leq k\leq n-1.$

So, we obtain
\[\beta_{x+e_k}=\beta_{x}, \quad \alpha_{x+e_k}=\alpha_{x}, \quad \beta_{e_2+e_k}=\beta_{e_k}, \quad 4\leq k\leq n-1, \quad \text{i.e.}\]
\[ \beta_{x}=\beta_{e_k}, \qquad 4\leq k\leq n-1.\]

Follows, the local automorphism $\Phi$ on $R_2$ has the next form:
 \[\left\{\begin{array}{ll}
     \Phi(e_1)=\alpha_{e_1} e_1+\sum\limits_{i=3}^n\frac{(-1)^i\alpha_{e_1}\beta_{e_1}^{i-2}}{(i-2)!}e_i, & \\
       \Phi(e_2)=\gamma_{e_2} e_2, & \\
      \Phi(e_i)=\sum\limits_{j=i}^n \frac{(-1)^{j-i}\alpha_{e_1}^{i-1} \beta_{e_1}^{j-i}}{(j-i)!} e_j,&  3\leq i\leq n,\\
       \Phi(x)=\beta_{e_1} e_1+\sum\limits_{i=3}^n\frac{(-1)^i\beta_{e_1}^{i-1}}{(i-1)!}e_i+x, & \\
        \Phi(y)=\delta_y e_2+y. & \\
     \end{array}\right.\]
By the definition of a local automorphism $\alpha_{e_1}\neq 0$ and $\beta_{e_2}\neq 0$.
which implies that $\Phi$ is an automorphism from theorem \ref{11}.
\end{proof}

\begin{theorem}
Every local automorphism on $R_3$ is an automorphism.
\end{theorem}

\begin{proof} The proof is similar to the proof of Theorem \ref{12}.
\end{proof}

	\section{2-local automorphisms of Solvable Leibniz algebras}\label{S:inf}

	\begin{theorem}
Every 2-local automorphism of $R_0$ is an automorphism.
\end{theorem}

\begin{proof}
\quad Let $\phi$ be an arbitrary 2 -local automorphism of $R_0$. Then, by the definition, for every element $x\in R_0$,
\[x=\sum\limits_{i=0}^n x_i e_i,\] \\
there exist elements $\alpha_{x,e_1}, \beta_{x,e_1}$ such that
\[A_{x,e_1}=\begin{pmatrix}
        1&0&0&\dots&0&0\\
        \alpha_{x,e_1}&\beta_{x,e_1}&0&\dots&0&0\\
        \frac{\alpha_{x,e_1}^2}{2}&\alpha_{x,e_1}\beta_{x,e_1}&\beta_{x,e_1}^2&\dots&0&0\\
        \vdots&\vdots&\vdots&\ddots&\vdots&\vdots\\
        \frac{\alpha_{x,e_1}^{n-1}}{(n-1)!}&\frac{\alpha_{x,e_1}^{n-2}\beta_{x,e_1}}{(n-2)!}&\frac{\alpha_{x,e_1}^{n-3}\beta_{x,e_1}^2}{(n-3)!}&\dots&\beta_{x,e_1}^{n-1}&0\\
        \frac{\alpha_{x,e_1}^n}{n!}&\frac{\alpha_{x,e_1}^{n-1}\beta_{x,e_1}}{(n-1)!}&\frac{\alpha_{x,e_1}^{n-2}\beta_{x,e_1}^2}{(n-2)!}&\dots&\alpha_{x,e_1}\beta_{x,e_1}^{n-1}&\beta_{x,e_1}^n
\end{pmatrix},\] \\
$\phi(x)=A_{x,e_1} \bar{x}$, where $\bar{x} = (x_0, x_1, x_2,\dots, x_n)^T$ is the vector corresponding
to $x$, and
\[\phi(e_1)=A_{x,e_1}\overline{e_1}=(
        0,\beta_{x,e_1},\alpha_{x,e_1}\beta_{x,e_1},\dots,\frac{\alpha_{x,e_1}^{n-2}\beta_{x,e_1}}{(n-2)!},\frac{\alpha_{x,e_1}^{n-1}\beta_{x,e_1}}{(n-1)!})^T.\]

Since $\phi(e_1) = \varphi_{x,e_1}(e_1) = \varphi_{y,e_1}(e_1)$, we have

\[\phi(e_1)=(0,\beta_{x,e_1},\alpha_{x,e_1}\beta_{x,e_1},\dots,\frac{\alpha_{x,e_1}^{n-2}\beta_{x,e_1}}{(n-2)!},\frac{\alpha_{x,e_1}^{n-1}\beta_{x,e_1}}{(n-1)!})^T=\]
\[=(0,\beta_{y,e_1},\alpha_{y,e_1}\beta_{y,e_1},\dots,\frac{\alpha_{y,e_1}^{n-2}\beta_{y,e_1}}{(n-2)!},\frac{\alpha_{y,e_1}^{n-1}\beta_{y,e_1}}{(n-1)!})^T\]
for each pair, $x,y$ of elements in $R_0$. Hense, $\alpha_{x,e_1}=\alpha_{y,e_1}, \beta_{x,e_1}=\beta_{y,e_1}$. Therefore
\[\phi(x)=A_{y,e_1}\bar{x}\] \\
for any $x\in R_0$, and the matrix of $\phi(x)$ does not depend on $x$. Thus, by theorem \ref{11}, $\phi$ is an automorphism.
\end{proof}

\begin{theorem}
Every 2-local automorphism of $R_1$ is an automorphism.
\end{theorem}
\begin{proof}

Let $z=\sum\limits_{i=1}^n z_i e_i +z_{n+1}x+z_{n+2}y $ be an arbitrary element from $ R_1$.
For every $v\in  R_1$ there exists an automorphism $\varphi_{v,z}$ such that
$$
\phi(v)=\varphi_{v,z}(v),\quad
\phi(z)=\varphi_{v,z}(z).
$$
Let $A_{v,z}=(a_{i,j}^{v,z})_{i,j=1}^{n+2}$
be the matrix of the automorphism $\varphi_{v,z}$.

Then from
$$
\varphi_{e_1,v}(e_1)=\varphi_{e_1,z}(e_1),\ v\in  R_1
$$
it follows that
$$
\alpha_{e_1,v}e_1=\alpha_{e_1,z}e_1,\ v\in  R_1.     \eqno{(2.1)}
$$
Hence, $\alpha_{e_1,v}=\alpha_{e_1,z}$. In particular, $\alpha_{e_1,e_2}=\alpha_{e_1,e_3}$.

Then from
$$
\varphi_{e_2,v}(e_2)=\varphi_{e_2,z}(e_2),\ v\in R_1
$$
it follows that
$$
\sum\limits_{i=2}^n \frac{(-1)^{i-2}\beta_{e_2,v}(\gamma_{e_2,v})^{i-2}}{(i-2)!}e_i=
\sum\limits_{i=2}^n \frac{(-1)^{i-2}\beta_{e_2,z}(\gamma_{e_2,z})^{i-2}}{(i-2)!}e_i
$$
Hence,
$$
\beta_{e_2,v}=\beta_{e_2,z}, \gamma_{e_2,v}=\gamma_{e_2,z}.
$$
In particular,
$$
\beta_{e_2,e_1}=\beta_{e_2,e_3}, \gamma_{e_2,e_1}=\gamma_{e_2,e_3}.
$$
For any $4\leq i\leq n-1$ we get
$$
\varphi_{e_{i-1},e_i}(e_i)=\varphi_{e_i,e_{i+1}}(e_i),
$$
and
$$
(\alpha_{e_{i-1},e_i})^{i-2}\beta_{e_{i-1},e_i}e_i=(\alpha_{e_i,e_{i+1}})^{i-2}\beta_{e_i,e_{i+1}}e_i,
$$
$$
-(\alpha_{e_{i-1},e_i})^{i-2}\beta_{e_{i-1},e_i}\gamma_{e_{i-1},e_i}e_{i+1}=-(\alpha_{e_i,e_{i+1}})^{i-2}\beta_{e_i,e_{i+1}}\gamma_{e_i,e_{i+1}}e_{i+1}.
$$
Hence, for any $4\leq i\leq n-1$, we get
$$
\gamma_{e_{i-1},e_i}=\gamma_{e_i,e_{i+1}}=\gamma_{e_1,x}
$$
by (2.1).
Also, by (2.1) we get
\[
\alpha_{e_1,e_i}=\alpha_{e_1,e_{i+1}}.
\]
Hence,
\[
\alpha_{e_{i-1},e_i}=\alpha_{e_1,e_i}=\alpha_{e_1,e_{i+1}}=\alpha_{e_i,e_{i+1}}
\]
and
\[
\beta_{e_{i-1},e_i}=\beta_{e_i,e_{i+1}}
\]
for any $4\leq i\leq n-1$.

Therefore, for every $i$ in $\{1,2,3,...,n-1\}$, the matrix $A_{e_i,e_{i+1}}=(a_{j,k})_{j,k=1}^{n+2}$
of the automorphism $\varphi_{e_i,e_{i+1}}$ is equal to the following matrix
\begin{tiny}
\[A=
\left(
  \begin{array}{ccccccccc}
\alpha_{e_1,e_2}&  0  &  0    & \cdots & 0 & \gamma_{e_1,x} & 0 \\
    0 & \beta_{e_1,e_2}  & 0  & \cdots & 0 &    0   & 0 \\
    0 &  -\beta_{e_1,e_2}\gamma_{e_1,x} & \alpha_{e_1,e_2}\beta_{e_1,e_2} & \cdots & 0 & 0 & 0 \\

    \vdots & \vdots & \vdots & \vdots & \vdots  &  \vdots & \vdots  \\
    0 & \frac{(-1)^{n-2}\beta_{e_1,e_2}(\gamma_{e_1,x})^{n-2}}{(n-2)!} &\frac{(-1)^{n-3}\alpha_{e_1,e_2}\beta_{e_1,e_2}(\gamma_{e_1,x})^2}{(n-3)!} & \cdots & (\alpha_{e_1,e_2})^n\beta_{e_1,e_2} &  0   &  0   \\
    0 & 0 & 0 & \cdots & 0 & 1 & 0 \\
    0 & 0 & 0 & \cdots & 0 & 0 & 1 \\
  \end{array}
\right)
\].
\end{tiny}

Let $v=e_1$, $z=e_1+e_2$. Then, from $\varphi_{e_1,e_2}(e_1)=\varphi_{e_1,e_1+e_2}(e_1)$ it follows that
\[
\alpha_{e_1,e_1+e_2}=\alpha_{e_1,e_2}.
\]
Now, let $v=e_2$, $z=e_1+e_2$. Then, from $\varphi_{e_1,e_2}(e_2)=\varphi_{e_2,e_1+e_2}(e_2)$ it follows that
\[
\beta_{e_2,e_1+e_2}=\beta_{e_1,e_2}, \beta_{e_2,e_1+e_2}\gamma_{e_2,e_1+e_2}=\beta_{e_1,e_2}\gamma_{e_1,x}, \gamma_{e_2,e_1+e_2}=\gamma_{e_1,x}.
\]
But
\[
\varphi_{e_1,e_1+e_2}(e_1+e_2)=\varphi_{e_2,e_1+e_2}(e_1+e_2)
\]
and
\[
\beta_{e_2,e_1+e_2}=\beta_{e_1,e_1+e_2}, \beta_{e_2,e_1+e_2}\gamma_{e_2,e_1+e_2}=\beta_{e_1,e_1+e_2}\gamma_{e_1,e_1+e_2}, \gamma_{e_2,e_1+e_2}=\gamma_{e_1,e_1+e_2}.
\]
Hence,
\[
\beta_{e_1,e_1+e_2}=\beta_{e_1,e_2}, \gamma_{e_1,e_1+e_2}=\gamma_{e_1,x}.
\]

Now, we take $v=e_1+e_2$. Then, from $\varphi_{e_1,e_1+e_2}(e_1+e_2)=\varphi_{e_1+e_2,z}(e_1+e_2)$ it follows that
\[
\alpha_{e_1,e_1+e_2}=\alpha_{e_1+e_2,z}.
\]
\[
\beta_{e_1,e_1+e_2}=\beta_{e_1+e_2,z}, \beta_{e_1,e_1+e_2}\gamma_{e_1,e_1+e_2}=\beta_{e_1+e_2,z}\gamma_{e_1+e_2,z}, \gamma_{e_1,e_1+e_2}=\gamma_{e_1+e_2,z}.
\]
Hence,
\[
\alpha_{e_1+e_2,z}=\alpha_{e_1,e_2}, \beta_{e_1+e_2,z}=\beta_{e_1,e_2}, \gamma_{e_1+e_2,z}=\gamma_{e_1,x}.
\]

So, the matrix of $\varphi_{e_1+e_2,z}$ coincides with the matrix $A$ for an arbitrary element $z$.
Note that $\alpha_{e_1,e_2}\neq 0$, $\beta_{e_1,e_2}\neq 0$ by the definition of a 2-local automorphism and theorem \ref{11}.
Hence, the 2-local automorphism $\phi$ is an automorphism. This ends the proof.
\end{proof}

\medskip

\begin{theorem} \label{13}
Every 2-local automorphism of $R_2$ is an automorphism.
\end{theorem}

\begin{proof}
Let $\phi$ be an arbitrary 2 -local automorphism of $R_2$. By the definition, for all $z,t\in R_2$,
there exists an automorphism $\varphi_{z,t} $ of $R_2$ such that
\[\phi(z)=\varphi_{z,t}(z), \quad \phi(t)=\varphi_{z,t}(t).\]
By theorem \ref{11}, the automorphism $\varphi_{z,t}$ has a matrix of the following form:
\[A_{z,t}=\begin{pmatrix}
        \alpha_{z,t}&0&0&0&\dots&0&\beta_{z,t}&0\\
        0&\gamma_{z,t}&0&0&\dots&0&0&\delta_{z,t}\\
        -\alpha_{z,t}\beta_{z,t}&0&\alpha_{z,t}^2&0&\dots&0&-\frac{\beta_{z,t}^2}{2}&0\\
         \frac{\alpha_{z,t}\beta_{z,t}^2}{2}&0&-\alpha_{z,t}^2\beta_{z,t}&\alpha_{z,t}^3&\dots&0&\frac{\beta_{z,t}^3}{6}&0\\
        \vdots&\vdots&\vdots&\vdots&\ddots&\vdots&\vdots&\vdots\\
         \frac{(-1)^n\alpha_{z,t}\beta_{z,t}^{n-2}}{(n-2)!}&0&\frac{(-1)^{n-3}\alpha_{z,t}^2\beta_{z,t}^{n-3}}{(n-3)!}&\frac{(-1)^{n-4}\alpha_{z,t}^3\beta_{z,t}^{n-4}}{(n-4)!}&\dots&\alpha_{z,t}^{n-1}&\frac{(-1)^n\beta_{z,t}^{n-1}}{(n-1)!}&0\\
         0&0&0&0&\dots&0&1&0\\
           0&0&0&0&\dots&0&0&1
\end{pmatrix}.\]

\quad In accordance with the equalities
\[\phi(e_1)=\varphi_{e_1,t}(e_1)=\varphi_{e_1,z}(e_1)\]
we obtain
\[(\alpha_{e_1,t},0,-\alpha_{e_1,t}\beta_{e_1,t},\frac{\alpha_{e_1,t}\beta_{e_1,t}^2}{2},\dots,\frac{(-1)^n\alpha_{e_1,t}\beta_{e_1,t}^{n-2}}{(n-2)!},0,0)^T=\]
\[=(\alpha_{e_1,z},0,-\alpha_{e_1,z}\beta_{e_1,z},\frac{\alpha_{e_1,z}\beta_{e_1,z}^2}{2},\dots,\frac{(-1)^n\alpha_{e_1,z}\beta_{e_1,z}^{n-2}}{(n-2)!}, 0,0)^T,\]
which implies
\[\alpha_{e_1,z}=\alpha_{e_1,t}, \beta_{e_1,z}=\beta_{e_1,t}.\]

\quad Considering the equality
\[\varphi_{e_2,z}(e_2)=\varphi_{e_2,t}(e_2)\]
we find that
\[\gamma_{e_2,z}=\gamma_{e_2, t}.\]

\quad Similarly, from
\[\varphi_{y,z}(y)=\varphi_{y,t}(y)\]
it follows that
\[\delta_{y,z}=\delta_{y,t}.\]

Hence,
\[\phi(z)=\varphi_{e_1,z}(z)=\varphi_{e_2,z}(z)=\varphi_{y,z}(z)\]
for any $z\in R_2$, and the matrix of $\phi(z)$ does not depend on $z$. Thus, by theorem \ref{11}, $\phi$ is an automorphism.
\end{proof}

\begin{theorem}
Every 2-local automorphism of $R_3$ is an automorphism.
\end{theorem}

\begin{proof}
 The proof is similar to the proof of Theorem \ref{13}.
\end{proof}

Summarizing and concluding the results on the paper we present the next conjecture:

\begin{conjecture}
Each local and 2-local automorphisms on the solvable Leibniz algebras with a given nilradical,
the dimension of whose complementary space is maximal, are automorphisms.
\end{conjecture}

\end{document}